\documentclass[oneside]{article}
\usepackage{amsmath,amssymb,amscd,amsthm,verbatim,alltt,amsfonts,array}
\usepackage{mathrsfs}
\usepackage[english]{babel}
\usepackage{latexsym}
\usepackage{amssymb}
\usepackage{euscript}
\usepackage{graphicx}
\usepackage{arydshln}
\usepackage[dvipsnames]{xcolor}

\newtheorem{theorem}{Theorem}
\theoremstyle{plain}

\newtheorem{remark}{Remark}

\numberwithin{equation}{section}

\begin{document}

{\footnotesize%
\hfill
}

  \vskip 1.2 true cm

\begin{center} {\bf  Stability in the homology of configuration spaces of high dimensional manifolds} \\
          {by}\\
{\sc Muhammad Yameen}
\end{center}

\pagestyle{myheadings}
\markboth{Stability in the homology of configuration spaces of high dimensional manifolds}{Muhammad Yameen}

\begin{abstract}
Consider the configuration spaces of manifold (closed or open). An influential theorem of McDuff and Segal shows that the (co)homology of unordered configuration spaces of open manifold is independent of number of configuration points in a range of degree called the stable range. Church prove that the homological stability is also present in the closed background manifold. The Church's stable range is only depend on the number of configuration points. In this paper, we prove that if the cohomology of manifold (not necessarily orientable) is concentrated in even degree than the stable range for odd degree homology groups of configuration spaces is also depend on the dimension of background manifold. Moreover, we prove that if the second cohomology group of manifold is nontrivial then our stable range for even degree homology groups is optimal, and this stable range is depend only on the number of configuration points.
\end{abstract}

\begin{quotation}

\end{quotation}

\thispagestyle{empty}

\section{Introduction}

\label{sec:intro}


For any manifold $M$, let
$$F_{k}(M):=\{(x_{1},\ldots,x_{k})\in M^{k}| x_{i}\neq x_{j}\,for\,i\neq j\}$$
be the configuration space of $k$ distinct ordered points in $M$ with induced topology. The symmetric group $S_{k}$ acts on $F_{k}(M)$ by permuting the coordinates. The quotient $$C_{k}(M):=F_{k}(M)/S_{k} $$
is the unordered configuration space with quotient topology. The geometry and topology of the configuration spaces have attracted an important attention, over the years. Configuration spaces is an active area of research of geometry with important applications in several areas including mathematical physics and computer science.

An important aspect of the study of configuration spaces of manifold is the phenomenon
of \textbf{homological stability}. Arnold proved that integral cohomology groups of configuration spaces of Euclidean plane are stabilized:
$$ H^{i}(C_{2i-2}(\mathbb{R}^{2});\mathbb{Z})\cong H^{i}(C_{2i-1}(\mathbb{R}^{2});\mathbb{Z})  \cong H^{i}(C_{2i}(\mathbb{R}^{2});\mathbb{Z})\cong\ldots $$
The isomorphisms (for $k$ large depending on $i$)
$$ H^{i}(C_{k}(M);\mathbb{Z})\cong H^{i}(C_{k+1}(M);\mathbb{Z})\cong H^{i}(C_{k+2}(M);\mathbb{Z})\cong\ldots $$
were generalized for open manifolds by D. McDuff \cite{MD} and G. Segal \cite{S}. The corresponding theorem is false for closed manifolds \cite{N}; for example $$H_{1}(C_{k}(S^{2});\mathbb{Z})=\mathbb{Z}_{2k-2}.$$ Using representation stability \cite{CF}, Church \cite{C} proved that
$$ H^{i}(C_{k}(M);\mathbb{Q})\cong H^{i}(C_{k+1}(M);\mathbb{Q})\cong H^{i}(C_{k+2}(M);\mathbb{Q})\cong\ldots $$
for $k>i$ and $M$ a connected oriented manifold of finite type. This result was extended by Randal-Williams \cite{RW}.
\begin{theorem}\label{Stable}
Let $M$ be a connected manifold of dimension bigger than 2. The cap product with unit in $H^{0}(M;\mathbb{Q})$ induces a map
$$ H_{i}(C_{k+1}(M);\mathbb{Q}) \rightarrow  H_{i}(C_{k}(M);\mathbb{Q})$$ 
that is an isomorphism for $i\leq k.$
\end{theorem} 
 We will study the stability properties of the rational (co)homology of unordered configuration spaces of connected even cohomology manifolds of higher dimension. We prove that if the cohomology of high dimensional manifold is concentrate in even degree then the stable range for odd degree homology groups is also depend on the dimension of manifold. We will restrict our attention to the case
of even dimensional manifolds, since there is a simple closed form available for odd dimensional manifolds \cite{BCT} (see last section for
further commentary). A manifold $M$ has {\em even cohomology} if all its odd cohomology groups are trivial.  Without a special mention, the (co)homology will have coefficients in $\mathbb{Q}.$
\begin{theorem}\label{stable}
Let $M$ be a connected even cohomology manifold of even dimension $d\geq6.$ The cap product with unit in $H^{0}(M;\mathbb{Q})$ induces a map
$$ H_{i}(C_{k+1}(M);\mathbb{Q}) \rightarrow  H_{i}(C_{k}(M);\mathbb{Q})$$ 
that is \\
$\bullet$ an isomorphism for $i\leq 2k$ when $i$ is even, and \\
$\bullet$ an isomorphism for $i\leq 2k+d-5$ when $i$ is odd.
\end{theorem} 
 There is no optimal stable range present in literature for configuration space of high dimensional manifolds. We prove that if the second cohomology group of manifold is non-trivial then our stable range for even degree homology group is optimal.
\begin{theorem}\label{optimal}
Let $M$ be a connected even cohomology  manifold of even dimension bigger than 5. If $M$ is not closed and $H^{2}(M;\mathbb{Q})$ is not trivial then the stable range for even degree homology group in Theorem \ref{stable} is optimal.
\end{theorem} 
\begin{remark}
\emph{ Our stable range can be further improved
if the low-degree Betti numbers of manifold vanish.}
\end{remark}


\section{Lie algebra homology}
The proof of Church is depend on the representation stability and Totaro's spectral sequence \cite{T}. However, we use the Chevalley–Eilenberg complex (for more details see \cite{F-M}\cite{K}\cite{F-Th} \cite{F-Ta}\cite{Kn}). \\\\
Let us introduce some notations. The symmetric algebra $Sym(W^{*})$ is the tensor product of a polynomial algebra and an exterior algebra:
$$ Sym(W^{*})=\bigoplus_{k\geq0}Sym^{k}(W^{*})=Poly(W^{even})\bigotimes Ext(W^{odd}), $$
where $Sym^{k}$ is generated by the monomials of length $k.$ The $n$-th suspension of the graded vector space $W$ is the graded vector space $W[n]$ with
$W[n]_{i} = W_{i-n},$ and the element of $W[n]$ corresponding to $a\in W$ is denoted $s^{n}a.$ We write $H_{-*}(M;\mathbb{Q})$ for the graded vector space whose degree $-i$ part is
the $i$-th homology group of $M.$

Now consider two graded vector spaces $$U^{*}=H^{-*}_{c}(M;\mathbb{Q}^{w})[d],\quad V^{*}=H_{c}^{-*}(M;\mathbb{Q})[2d-1]:$$
where
$$ U^{*}=\bigoplus_{i=0}^{d}U^{i},\quad V^{*}=\bigoplus_{j=d-1}^{2d-1}V^{j}.$$
We choose bases in $U^{i}$ and $V^{j}$ as 
$$U^{i}=\mathbb{Q}\langle u_{i,1},u_{i,2},\ldots\rangle,\quad V^{j}=\mathbb{Q}\langle v_{j,1},v_{j,2},\ldots\rangle$$
(the degree of an element is marked by the first lower index; $x_{i}^{l}$ stands for the product $x_{i}\wedge x_{i}\wedge\ldots\wedge x_{i}$ of $l$-factors). Always we take $U^{0}=\mathbb{Q}\langle u_{0}\rangle$. Now consider the graded algebra
$$ \Theta^{*,*}_{k}(M)=\bigoplus_{i\geq 0}\bigoplus_{\omega=0}^{\left\lfloor\frac{k}{2}\right\rfloor}
\Theta^{i,\omega}_{k}(M)=\bigoplus_{\omega=0}^{\left\lfloor\frac{k}{2}\right\rfloor}\,(Sym^{k-2\omega}(U^{*})\otimes Sym^{\omega}(V^{*})) $$
where the total degree $i$ is given by the grading of $U^{*}$ and $V^{*}$. We called $\omega$ is a weight grading. The differential $\partial:Sym^{2}(U^{*})\rightarrow V^{*}$ is defined as a coderivation by the equation 
$$\partial(s^{d}a\wedge s^{d}b)=(-1)^{(d-1)|b|}s^{2d-1}(a\cup b),$$ where $$\cup\,:H^{-*}_{c}(M;\mathbb{Q}^{w})^{\otimes2}\rightarrow H^{-*}_{c}(M;\mathbb{Q})$$
(here $H^{-*}_{c}$ denotes compactly supported cohomology of $M$). The degree of $\partial$ is $-1.$ The differential $\partial$ extends over  $\Theta^{*,*}_{k}(M)$ by co-Leibniz rule. By definition the elements in $U^{*}$ have length 1 and weight 0 and the elements in $V^{*}$ have length 2 and weight 1. By definition of differential, we have 
$$\partial:\Theta^{*,*}_{k}(M)\longrightarrow\Theta^{*-1,*+1}_{k}(M).$$

\begin{theorem}
	If $d$ is even, $H_{*}(C_{k}(M);\mathbb{Q})$  is isomorphic to the homology of the complex 
	$$ (\Theta^{*,*}_{k}(M),\partial).$$
\end{theorem}

\section{Proofs of the main results}
In this section, we give the proofs of main theorems.\\\\
\textit{Proof of Theorem \ref{stable}.} Let $d$ be even and $M$ is even cohomology manifold.\\\\
\textbf{Case I.} Assume $M$ is is not closed. We have decomposition of algebras:
$$\Theta^{*,*}_{k+1}(M)\cong (\langle u_{0}\rangle\otimes \Theta^{*,*}_{k}(M))\oplus \Theta^{*,*}_{k+1}(M)/\langle u_{0}\rangle.$$
The subspace $\Theta^{*,*}_{k+1}(M)/\langle u_{0}\rangle$ is $\partial-$invariant:
$$\partial(\Theta^{*,*}_{k+1}(M)/\langle u_{0}\rangle)\subseteq \Theta^{*,*}_{k+1}(M)/\langle u_{0}\rangle.$$
The degrees of elements of $U^{*}$ are less than $d.$ Also, the degree of elements of $V^{*}$ are greater than and equal to $d-1.$ The differential $\partial$ has degree $-1.$ Therefore, $\partial(\langle u_{0}\rangle\otimes U^{*})=0.$ 
This implies that the subspace $\langle u_{0}\rangle\otimes \Theta^{*,*}_{k}(M)$ is also $\partial-$invariant:
$$\partial(\langle u_{0}\rangle\otimes \Theta^{*,*}_{k}(M))\subseteq \langle u_{0}\rangle\otimes \Theta^{*,*}_{k}(M).$$
As a consequence, we get the decomposition of complexes:
$$(\Theta^{*,*}_{k+1}(M),\partial)\cong (\langle u_{0}\rangle\otimes \Theta^{*,*}_{k}(M),\partial)\oplus (\Theta^{*,*}_{k+1}(M)/\langle u_{0}\rangle,\partial).$$
Moreover, we have an isomorphism :
$$H_{*}(C_{k}(M);\mathbb{Q})\cong (\langle u_{0}\rangle\otimes \Theta^{*,*}_{k}(M),\partial)$$
Hence, we get an isomorphism :
$$H_{*}(C_{k+1}(M);\mathbb{Q})\cong H_{*}(C_{k}(M);\mathbb{Q})\oplus H_{*}(\Theta^{*,*}_{k+1}(M)/\langle u_{0}\rangle,\partial).$$
Moreover, for each non-zero $u\in \Theta^{*,*}_{k+1}(M)/\langle u_{0}\rangle,$ we have
$$ \mbox{deg(u)} \geq\begin{cases}
      2k+2, & \mbox{ if deg(u) is even}\\
      2k+d-3, & \mbox{ if deg(u) is odd}.\\
   \end{cases}
$$
As a consequence, we conclude that the map
$$ H_{i}(C_{k+1}(M);\mathbb{Q}) \rightarrow  H_{i}(C_{k}(M);\mathbb{Q})$$ 
 is an isomorphism for $i\leq 2k$ when $i$ is even, and an isomorphism for $i\leq 2k+d-5$ when $i$ is odd.\\\\
\textbf{Case II.} Assume $M$ is closed. The homology of $C_{k}(M)$ is stabilized in the stable range of $H_{*}(C_{k}(M-\{pt\});\mathbb{Q})$ (see Theorem 3 of \cite{RW}). This completes the proof.
$\hfill \square$\\\\
\textit{Proof of Theorem \ref{optimal}.}
Let $d\geq 6$ be even and $M$ is even cohomology manifold (but not closed). Let $\mbox{dim} H^{2}(M;\mathbb{Q})=n.$ Then the subspace $U^{2}$ of $U^{*}$ is 
$$U^{2}=\langle u_{2,1},\ldots,u_{2,n}\rangle.$$ The degree of elements of $V^{*}$ are greater than 4. The differential has bi-degree $(-1,1).$ Therefore, $\partial(u_{2,i}\otimes u_{2,j})=0.$ Each $u_{2,i_{1}}\ldots u_{2,i_{k+1}}\in \text{I}^{k+1}\subseteq \Theta^{2k+2,0}/\langle u_{0}\rangle$ gives a cohomology  class, where $\text{I}^{k+1}=\langle u_{2,i_{1}}\ldots u_{2,i_{k+1}}|\, i_{j}\in\{1,\ldots n\}\rangle.$ The carnality of the bases elements of $\text{I}^{k+1}$ is $\binom{(k+1)n-1}{n-1}.$ From the proof of Theorem \ref{Stable}, we have decomposition 
$$H_{*}(C_{k+1}(M);\mathbb{Q})\cong H_{*}(C_{k}(M);\mathbb{Q})\oplus H_{*}(\Theta^{*,*}_{k+1}(M)/\langle u_{0}\rangle,\partial).$$
This implies that for each $i\geq0,$ we have
$$ \mbox{dim} H_{i}(C_{k+1}(M);\mathbb{Q})\geq  \mbox{dim} H_{i}(C_{k}(M);\mathbb{Q}).$$ 
In particular, for $i=2k+2$ we have 
$$ \mbox{dim} H_{2k+2}(C_{k+1}(M);\mathbb{Q})\geq \mbox{dim} H_{2k+2}(C_{k}(M);\mathbb{Q})+\binom{(k+1)n-1}{n-1}.$$ 
This implies that for $n\neq0$ we have
$$ \mbox{dim} H_{2k+2}(C_{k+1}(M);\mathbb{Q})> \mbox{dim} H_{2k+2}(C_{k}(M);\mathbb{Q}).$$
As a consequence, we conclude that the homology groups $H_{2k+2}(C_{k+1}(M);\mathbb{Q})$ and  $H_{2k+2}(C_{k}(M)$ are not isomorphic.
$\hfill \square$\\\\
\section{Odd dimensional case}
We comment briefly on the case of odd dimensional manifolds, which is encompassed
implicitly by the results of \cite{BCT}.  We have an isomorphism \cite{BCT}:
$$H_{*}(C_{k}(M);\mathbb{Q})\cong Sym^{k}(H_{*}(M;\mathbb{Q})).$$ 
We can write the decomposition of algebras: 
$$Sym^{k}(U^{*})\cong (\langle u_{0}\rangle\otimes Sym^{k-1}(U^{*}))\oplus Sym^{k}(U^{*})/\langle u_{0}\rangle,$$
here we take $U^{*}\cong H_{*}(M;\mathbb{Q})$ and $u_{0}$ is a corresponding element of $H_{0}(M;\mathbb{Q}).$
Moreover, we have an isomorphism :
$$H_{*}(C_{k-1}(M);\mathbb{Q})\cong \langle u_{0}\rangle\otimes Sym^{k-1}(U^{*}).$$
As a consequence, we get 
$$H_{*}(C_{k}(M);\mathbb{Q})\cong H_{*}(C_{k-1}(M);\mathbb{Q})\oplus Sym^{k}(U^{*})/\langle u_{0}\rangle.$$
This implies that for each $i\geq0,$ the map
$$k\rightarrow \emph{dim} H_{i}(C_{k}(M);\mathbb{Q})$$ 
is eventually constant. If $M$ has even cohomology than the odd degree homology groups of $C_{k}(M)$ are vanish.\\\\

\noindent\textbf{Acknowledgement}\textit{.} The author gratefully acknowledge the support from the ASSMS, GC university Lahore. This research is partially supported by Higher Education Commission of Pakistan.

\vskip 0,65 true cm




\vskip 0,65 true cm







\null\hfill  Abdus Salam School of Mathematical Sciences,\\
\null\hfill  GC University Lahore, Pakistan. \\
\null\hfill E-mail: {yameen99khan@gmail.com}

\end{document}